\documentclass[11pt]{article}

\usepackage{oldlfont}
\usepackage{amsfonts}
\usepackage{graphicx}
\usepackage{eepic}

\newcommand{\be}{\begin{equation}}
\newcommand{\ee}{\end{equation}}
\newcommand{\bea}{\begin{eqnarray}}
\newcommand{\eea}{\end{eqnarray}}
\newcommand{\beas}{\begin{eqnarray*}}
\newcommand{\eeas}{\end{eqnarray*}}
\newcommand{\ba}{\begin{array}}
\newcommand{\ea}{\end{array}}

\newcommand{\<}  {\langle}
\renewcommand{\>}{\rangle}
\newcommand{\field}[1]{\mathbb{#1}}

\newcommand{\diam}{{\mathrm diam}}
\newcommand{\grad}{\nabla}
\renewcommand{\div}{{\mathrm div}}
\newcommand{\curl}{{\mathrm curl}}

\newcommand{\G}{\Gamma}
\newcommand{\tG}{\tilde\Gamma}
\newcommand{\g}{\gamma}

\newcommand{\bn}{{\mathbf n}}
\newcommand{\tbn}{{\tilde{\mathbf n}}}

\newcommand{\bu}{{\mathbf u}}
\newcommand{\bv}{{\mathbf v}}
\newcommand{\bw}{{\mathbf w}}
\newcommand{\bx}{{\mathbf x}}

\newcommand{\boldf}{{\mathbf f}}

\newcommand{\bX}{{\mathbf X}}
\newcommand{\bH}{{\mathbf H}}
\newcommand{\bL}{{\mathbf L}}
\newcommand{\bcurl}{{\mathbf curl}}
\newcommand{\bzero}{{\mathbf 0}}

\newcommand{\gradg}{\nabla_{\G}}
\newcommand{\divg}{{\mathrm div}_{\G}}

\newcommand{\bcurlg}{{\mathbf curl}_{\G}}

\newcommand{\bnu}{\hbox{\mathversion{bold}$\nu$}}
\newcommand{\bxi}{\hbox{\mathversion{bold}$\xi$}}
\newcommand{\bPsi}{\hbox{\mathversion{bold}$\Psi$}}

\newcommand{\bphi}{\hbox{\mathversion{bold}$\phi$}}
\newcommand{\btau}{\hbox{\mathversion{bold}$\tau$}}

\newcommand{\CI}{{\cal I}}

\newcommand{\CM}{{\cal M}}
\newcommand{\CP}{{\cal P}}
\newcommand{\bCP}{\hbox{\mathversion{bold}$\cal P$}}
\newcommand{\CQ}{{\cal Q}}
\newcommand{\CT}{{\cal T}}
\newcommand{\CDN}{{\cal DN}}

\newcommand{\tH}{\tilde H}

\newtheorem{theorem}{Theorem} [section]
\newtheorem{lemma}{Lemma} [section]

\newtheorem{remark}{Remark} [section]
\newenvironment{proof}{\noindent\textbf{Proof.}\ }
              {\nopagebreak\hbox{ }\hfill$\Box$\bigskip}

\title{The $hp$-BEM with quasi-uniform meshes
       for the electric field integral equation on polyhedral surfaces:
       a priori error analysis
\thanks{Supported by EPSRC under grant no. EP/E058094/1.}}

\author{Alexei Bespalov
\thanks{Department of Mathematical Sciences, Brunel University,
        Uxbridge, West London UB8 3PH, UK.
        Email: {\tt albespalov@yahoo.com}
        }
        \and
        Norbert Heuer
\thanks{Facultad de Matem\'aticas, Pontificia Universidad Cat\'olica de Chile,
        Avenida Vicu\~na Mackenna 4860, Santiago, Chile.
        Email: {\tt nheuer@mat.puc.cl}}
        }
\begin{document}
\date{}
\maketitle

\begin{abstract}
This paper presents an a priori error analysis of the $hp$-version of the
boundary element method for the electric field integral equation on
a piecewise plane (open or closed) Lipschitz surface.
We use $\bH(\div)$-conforming discretisations with Raviart-Thomas elements
on a sequence of quasi-uniform meshes of triangles and/or parallelograms.
Assuming the regularity of the solution to the electric
field integral equation in terms of Sobolev spaces of tangential vector fields,
we prove an a priori error estimate of the method in the energy norm.
This estimate proves the expected rate of convergence with respect to
the mesh parameter $h$ and the polynomial degree $p$.

\bigskip
\noindent
{\em Key words}: $hp$-version with quasi-uniform meshes, boundary element method,
                 electric field integral equation,
                 time-harmonic electro-magnetic scattering,
                 a priori error estimate

\noindent
{\em AMS Subject Classification}: 65N38, 65N15, 78M15, 41A10
\end{abstract}

\section{Introduction} \label{sec_intro}
\setcounter{equation}{0}

With this paper we continue the analysis of high-order boundary element
methods (BEM) for the electric field integral equation (EFIE) started in
\cite{BespalovH_NpB,BespalovH_Chp}.
Our BEM is based on discretisations of the variational formulation
of the EFIE (called Rumsey's principle) with an $\bH(\div)$-conforming family
of boundary elements. This approach is referred to as the natural boundary element
method for the EFIE.
In \cite{BespalovH_NpB} we analysed the natural $p$-BEM for the EFIE on a
plane open surface with polygonal boundary. We proved convergence of the $p$-version
with Raviart-Thomas (RT) parallelogram elements
and derived an a priori error estimate which takes into account the strong
singular behaviour of the solution at edges and corners of the surface.
In our previous paper \cite{BespalovH_Chp} we considered the EFIE on
a piecewise plane (open or closed) Lipschitz surface $\G$ and proved
quasi-optimal convergence
of the natural $hp$-BEM with quasi-uniform meshes of triangles
and quadrilaterals. In the present note
we perform an a priori error analysis of that method on affine meshes
under the assumption that the regularity of the exact solution is given in
Sobolev spaces of tangential vector fields on $\G$. As the main result we prove
an a priori error estimate in the energy norm (Theorem~\ref{thm_converge}).
The estimate appears to be optimal with respect to the mesh size $h$ and the
polynomial degree $p$ as the convergence rates in both $h$ and $p^{-1}$ are $r+1/2$
for $p$ large enough. This corresponds to the expected rate which is the Sobolev
regularity order $r$ of the exact solution minus the Sobolev order
$-1/2$ of the energy norm.

While in the $h$-version the degrees of approximating polynomials
are fixed (usually at a low level) and convergence is achieved by refining the mesh,
the $p$-version keeps the mesh fixed and improves approximations
by increasing polynomial degrees. The $hp$-version combines both
mesh refinement and increase of polynomial degrees.
For boundary integral equations governing the Laplace equation optimal
$hp$-convergence rates for singular
problems (and quasi-uniform meshes) are proved in
\cite{BespalovH_08_hpB,BespalovH_hpW}. The analysis
of optimal $hp$-BEM convergence rates for the EFIE with singular solutions
is an open problem and under investigation. In this paper we deal with
the case of solutions with Sobolev regularity. We also note that the $hp$-BEM with
geometrically graded meshes (yielding an exponential rate of convergence)
has been studied in \cite{HeuerMS_99_ECB}, again for the Laplace equation
and for hypersingular and weakly singular integral operators. For the EFIE
its analysis is an open problem.

An a priori error analysis of the natural $h$-BEM for the EFIE was performed
in \cite{HiptmairS_02_NBE} for polyhedral surfaces and in \cite{BuffaC_03_EFI}
for open Lipschitz surfaces (see also \cite{BuffaH_03_GBE} for a survey of results
and techniques). In particular, an optimal $h$-convergence rate of the method for
given Sobolev regularity of the solution and given polynomial degree has been
proved in \cite{BuffaC_03_EFI}
(being proved on open surfaces, this result extends to polyhedral surfaces
as well). If the BEM for the EFIE converges quasi-optimally, then a priori error
analysis reduces to an approximation problem within the energy space,
which is either $\bH^{-1/2}(\divg,\G)$ or $\tilde\bH^{-1/2}_0(\divg,\G)$ depending on
whether the surface $\G$ is closed or open.
The main tool in the approximation analysis is an appropriate
$\bH(\div)$-conforming interpolation operator.
Whereas the standard $\bH(\div)$-conforming RT interpolation operator
works well for the $h$-version (cf. \cite{BuffaC_03_EFI}),
it does not provide an optimal result for high order methods.
The main reason for this is the lack of
stability (with respect to polynomial degrees) of this operator
for low-regular vector fields, which always appear when
dealing with the EFIE on non-smooth surfaces.
Furthermore, existing techniques to prove $p$-estimates for the
error of RT interpolation work only on quadrilateral elements,
and their extension to triangular elements does not seem feasible.

An alternative to the classical RT interpolation operator is a corresponding
projection based interpolation operator. Such operators were first introduced
in \cite{DemkowiczB_03_pIE} to analyse high-order finite element approximations
with $\bH(\curl)$-conforming edge elements for Maxwell's equations in two
dimensions. However, applying a simple rotation argument the results of
\cite{DemkowiczB_03_pIE} can be formulated in the $\bH(\div)$-conforming setting,
which is intrinsic to natural boundary element discretisations of the EFIE
(see also \cite{BespalovH_OEE}).
The projection based interpolation operators are stable with respect to
polynomial degrees, they work equally well on both
triangular and quadrilateral elements and also for low-regular fields.
That is why these operators have become an efficient tool in
the analysis of high-order methods (see \cite{BoffiDC_03_DCp,BoffiCDD_06_Dhp,Hiptmair_DCp}
for the finite element methods and \cite{BespalovH_NpB,BespalovH_Chp} for the BEM).
The a priori error analysis in this paper relies on such an operator as well.
In particular, we demonstrate that employing the $\bH(\div)$-conforming
projection based interpolation operator (rather than the
classical RT interpolation operator) one obtains an optimal error estimate
for the $hp$-BEM with quasi-uniform meshes.

The paper is organised as follows.
In the next section we formulate the EFIE (in a variational form) and
define the $hp$-version of the BEM with quasi-uniform meshes.
We also formulate the main result (Theorem~\ref{thm_converge}),
which states an a priori error estimate of the approximation method.
Section~\ref{sec_prelim} gives necessary preliminaries: first, in \S\ref{sec_spaces}
we introduce the needed notation and recall definitions of Sobolev spaces of
tangential vector fields;
then, in \S\ref{sec_int} we sketch the definition of the
$\bH(\div)$-conforming projection based interpolation operator on the reference
element and prove a new property of this operator related to
approximations of normal traces on the element's edges (Lemma~\ref{lm_aux1}).
In Section~\ref{sec_approx} we study approximating properties of the discrete
(boundary element) space $\bX_{hp}$ in the energy space $\bX$ of the EFIE, and
prove that the orthogonal projection onto $\bX_{hp}$ with respect to the norm
in $\bX$ satisfies an optimal error estimate in both $h$ and $p$.

Throughout the paper, $C$ denotes a generic positive constant which is independent
of $h$, $p$ and involved functions.

\section{Formulation of the problem and the main result} \label{sec_result}
\setcounter{equation}{0}

We consider the EFIE for which we have proved quasi-optimal convergence
of the $hp$-BEM with quasi-uniform meshes \cite{BespalovH_Chp}.
In this paper we provide an a priori error estimate. To this end let us recall
the model problem, its $hp$-discretisation and the involved spaces.

Let $\G$ denote a piecewise plane (open or closed) Lipschitz surface in ${\field{R}}^3$.
In the case of an open surface we additionally assume that $\G$ is orientable.
Let us introduce Rumsey's formulation of the electric field integral equation on $\G$.
For a given wave number $k>0$ and a scalar function $v$ (resp., tangential vector field $\bv$)
we define the single layer operator $\Psi_k$ (resp., $\bPsi_k$) by
\beas
     & \displaystyle{
     \Psi_k v(x) =
     \frac 1{4\pi}\int_\G v(y) \frac {e^{ik|x-y|}}{|x-y|}\,dS_y,\qquad
     x \in {\field{R}}^3 \backslash \G}
     &
     \\[5pt]
     & \displaystyle{
     \Big(\hbox{resp.,}\quad
     \bPsi_k \bv(x) =
     \frac 1{4\pi}\int_\G \bv(y) \frac {e^{ik|x-y|}}{|x-y|}\,dS_y,\qquad
     x \in {\field{R}}^3 \backslash \G\Big).}
     &
\eeas
Let $\bL^2_t(\G)$ be the space of two-dimensional, tangential, square integrable
vector fields on $\G$. By $\gradg$ (resp., $\divg$) we denote the surface gradient
(resp., surface divergence) acting on scalar functions
(resp., tangential vector fields) on $\G$. We will need the following space:
\[
  \bX = \bH^{-1/2}(\divg,\G) :=
        \{\bu \in \bH^{-1/2}_{\|}(\G);\; \divg\,\bu \in H^{-1/2}(\G)\}
\]
if $\G$ is a closed surface, and
\beas
     \bX = \tilde\bH^{-1/2}_0(\divg,\G)
     & := &
     \{\bu \in \tilde\bH^{-1/2}_{\|}(\G);\ \divg\,\bu \in \tH^{-1/2}(\G)\ \hbox{and}
     \\[3pt]
     &   &
     \quad
     \<\bu,\gradg v\> + \<\divg\,\bu,v\> = 0\quad
     \hbox{for all}\ \ v \in C^{\infty}(\G)\}
\eeas
if $\G$ is an open surface.
In the latter definition the brackets $\<\cdot,\cdot\>$ denote
dualities associated with $\bH^{1/2}_{\|}(\G)$ and $H^{1/2}(\G)$, respectively.
For definitions of the space $C^{\infty}(\G)$ and the Sobolev spaces on $\G$
we refer to \S\ref{sec_spaces} below.
Throughout, we use boldface symbols for vector fields.
The spaces (or sets) of vector fields are also denoted in boldface
(e.g., $\bH^s(\G) = (H^s(\G))^3$).

Let $\bX'$ be the dual space of $\bX$ (with $\bL^2_t(\G)$ as pivot space).
Now, for a given tangential vector field $\boldf\in \bX'$
($\boldf$ represents the excitation by an incident wave), Rumsey's formulation
reads as:
{\em find a complex tangential field $\bu\in\bX$ such that}
\be \label{bie_var}
    a(\bu,\bv)
    :=
    \<\g_{\rm tr}(\Psi_k\divg\,\bu),\divg\,\bv\> - k^2 \<\pi_{\tau}(\bPsi_k\bu), \bv\>
    = \<\boldf, \bv\> \quad\forall \bv \in \bX.
\ee
Here $\g_{\rm tr}$ is the standard trace operator, and $\pi_{\tau}$ denotes the
tangential components trace mapping (see \S\ref{sec_spaces} for the definition).
To ensure the uniqueness of the solution to
(\ref{bie_var}) in the case of the closed surface $\G$ we always assume
that $k^2$ is not an electrical eigenvalue of the interior problem.

For the approximate solution of (\ref{bie_var}) we apply the $hp$-version of the BEM based
on Galerkin discretisations with Raviart-Thomas spaces on quasi-uniform meshes.
In what follows, $h > 0$ and $p \ge 1$ will always specify the mesh parameter and
a polynomial degree, respectively.
For any $\Omega\subset {\field{R}}^n$ we will denote
$\rho_{\Omega} = \sup \{\diam (B);\; \hbox{$B$ is a ball in $\Omega$}\}$.
Furthermore, throughout the paper, $K$ is either the equilateral reference triangle
$T = \{x_2 > 0,\ x_2 < x_1 \sqrt 3,\ x_2 < (1-x_1)\sqrt 3\}$ or
the reference square $Q=(0,1)^2$. A generic side of $K$ will be denoted by $\ell$.

Let $\CT = \{\Delta_h\}$ be a family of meshes
$\Delta_h = \{\G_j;\; j=1,\ldots,J\}$ on $\G$, where the elements $\G_j$ are open
triangles or parallelograms such that $\bar\G = \cup_{j=1}^{J} \bar\G_j$,
and the intersection of any two elements $\bar\G_j,\,\bar\G_k$ ($j \not= k$)
is either a common vertex, an entire side, or empty.

We denote $h_j=\diam (\G_j)$ for any $\G_j\in \Delta_h$.
The elements are assumed to be shape regular, i.e., there exists a positive constant $C$
independent of $h = \max\limits_j h_j$ such that for
any $\G_j\in\Delta_h$ and arbitrary $\Delta_h\in \CT$ there holds
$h_j \le C\,\rho_{\G_j}$. Furthermore, any element $\G_j$ is the image
of the corresponding reference element $K$ under an affine mapping $T_j$,
more precisely
\[
  \bar\G_j = T_j(\bar K),\quad
  \bx = T_j(\bxi),\ 
  \bx = (x_1,x_2) \in \bar\G_j,\ \bxi = (\xi_1,\xi_2) \in \bar K.
\]
The Jacobian matrix of $T_j$ is denoted by $DT_j$ and its determinant
$J_j := \hbox{det}(DT_j)$ satisfies the relation $|J_j| \simeq h_j^2$.

We consider a family $\CT$ of quasi-uniform meshes $\Delta_h$
on $\G$ in the sense that there exists a positive constant
$C$ independent of $h$ such that for any $\G_j\in\Delta_h$ and
arbitrary $\Delta_h\in \CT$ there holds $h \le C\,h_j$.

The mapping $T_j$ introduced above is used to associate the scalar function $u$ defined
on the real element $\G_j$ with the function $\hat u$ defined on the reference
element $K$:
\[
  u = \hat u \circ T_j^{-1}\ \ \hbox{on $\G_j$\qquad and}\qquad
  \hat u = u \circ T_j\ \ \hbox{on $K$}.
\]
Any vector-valued function $\hat\bv$ defined on $K$ is transformed to the function
$\bv$ on $\G_j$ by using the Piola transformation:
\be \label{Piola}
    \bv = \CM_j(\hat \bv) = \hbox{$\frac {1}{J_j}$} DT_j \hat\bv \circ T_j^{-1},
    \quad
    \hat\bv = \CM_j^{-1}(\bv) = J_j DT_j^{-1} \bv \circ T_j.
\ee

Let us introduce the needed polynomial sets.
By $\CP_p(I)$ we denote the set of polynomials of degree $\le p$ on
an interval $I\subset\field{R}$, and $\CP^0_p(I)$ denotes the subset of $\CP_p(I)$
which consists of polynomials vanishing at the end points of $I$.
In particular, these two sets will be used for an edge $\ell \subset \partial K$.

Further, $\CP^1_p(T)$ denotes the set of polynomials on $T$ of total degree $\le p$, and
$\CP^2_{p_1,p_2}(Q)$ is the set of polynomials on $Q$ of degree $\le p_1$ in $\xi_1$
and degree $\le p_2$ in $\xi_2$. For $p_1 = p_2 = p$ we denote
$\CP^2_{p}(Q) = \CP^2_{p,p}(Q)$, and we will use the unified notation
$\CP_{p}(K)$, which refers to $\CP^1_{p}(T)$ if $K=T$ and to $\CP^2_{p}(Q)$ if $K=Q$.
The corresponding set of polynomial (scalar) bubble functions on $K$ is denoted
by $\CP^0_{p}(K)$.

Let us denote by $\bCP^{\rm RT}_p(K)$ the RT-space of order $p\ge 1$ on
the reference element $K$ (see, e.g., \cite{BrezziF_91_MHF, RobertsT_91_MHM}), i.e.,
\[
  \bCP^{\rm RT}_p(K) =
  (\CP_{p-1}(K))^2 \oplus \bxi \CP_{p-1}(K) =
  \cases{
         (\CP^1_{p-1}(T))^2 \oplus \bxi \CP^1_{p-1}(T)
         & \hbox{if \ $K = T$},
         \cr
         \noalign{\vspace{5pt}}
         \CP^2_{p,p-1}(Q) \times \CP^2_{p-1,p}(Q)
         & \hbox{if \ $K = Q$}.
         \cr
        }
\]
The subset of $\bCP^{\rm RT}_p(K)$ which consists of vector-valued polynomials
with vanishing normal trace on the boundary $\partial K$
(vector bubble-functions) will be denoted by $\bCP^{\rm RT,0}_p(K)$.

Then using transformations (\ref{Piola}), we set
\be \label{Xp}
    \bX_{hp} := \{\bv \in \bX^0;\;
                      \CM_j^{-1}(\bv|_{\G_j}) \in \bCP^{\rm RT}_p(K),\ j=1,\ldots,J\},
\ee
where the space $\bX^0 \subset \bX$ is defined in \S\ref{sec_spaces}
($\bX^0 = \bH(\divg,\G)$ if $\G$ is closed and
$\bX^0 = \bH_0(\divg,\G)$ if $\G$ is an open surface).
We will denote by $N = N(h,p)$ the dimension of the discrete space $\bX_{hp}$.
One has $N \simeq h^{-2}$ for fixed $p$ and $N \simeq p^2$ for fixed $h$.

The $hp$-version of the Galerkin BEM for the EFIE reads as:
{\em Find $\bu_{hp}\in \bX_{hp}$ such that}
\be \label{BEM}
    a(\bu_{hp},\bv) = \<\boldf, \bv\> \quad \forall \bv \in \bX_{hp}.
\ee

First, let us formulate the result which states the unique solvability of (\ref{BEM})
and quasi-optimal convergence of the $hp$-version of the BEM for the EFIE.

\begin{theorem} \label{thm_solve}
{\rm \cite[Theorem~2.1]{BespalovH_Chp}}
There exists $N_0 \ge 1$ such that for any $\boldf \in \bX'$ and for arbitrary
mesh-degree combination satisfying $N(h,p) \ge N_0$ the discrete problem {\rm (\ref{BEM})}
is uniquely solvable and the $hp$-version of the Galerkin BEM generated by RT-elements
converges quasi-optimally, i.e.,
\be \label{quasi-optimality}
    \|\bu - \bu_{hp}\|_{\bX} \le
    C \inf\{\|\bu - \bv\|_{\bX};\; \bv\in \bX_{hp}\}.
\ee
Here, $\bu \in \bX$ is the solution of {\rm (\ref{bie_var})},
$\bu_{hp} \in \bX_{hp}$ is the solution of {\rm (\ref{BEM})},
$\|\cdot\|_{\bX}$ denotes the norm in $\bX$, and
$C>0$ is a constant independent of $h$ and $p$.
\end{theorem}

The following theorem is the main result of this paper
(its formulation involves the space $\bX^r$ and the norm
$\|\cdot\|_{\bX^r}$ which are defined in~\S\ref{sec_spaces}).

\begin{theorem} \label{thm_converge}
Let $\bu \in \bX$ and $\bu_{hp} \in \bX_{hp}$ be the solutions
of {\rm (\ref{bie_var})} and {\rm (\ref{BEM})}, respectively.
Then there exists a real number $r$ ($r > 0$ if $\G$ is closed,
and $r \in (-\frac 12,0)$ if $\G$ is an open surface)
such that $\bu \in \bX^r$ and the following a priori error estimate holds
\be \label{convergence}
    \|\bu - \bu_{hp}\|_{\bX} \le
    C\,h^{1/2+\min\,\{r,p\}}\,p^{-(r+1/2)}\,\|\bu\|_{\bX^r}
\ee
with a positive constant $C$ independent of $h$ and $p$.
\end{theorem}

\begin{proof}
The assertion regarding the regularity of the solution $\bu$ to the EFIE is a
direct consequence of the regularity results in
\cite[Section 4.4]{CostabelD_00_SEF} (see also \cite[Appendix~A]{BespalovH_NpB}).
Due to the quasi-optimal convergence (\ref{quasi-optimality}) of the $hp$-BEM,
the error estimate in (\ref{convergence}) then immediately follows from the approximation
result of Theorem~\ref{thm_main_approx} below.
\end{proof}

\begin{remark} \label{rem_reg}
If $\G$ is a closed Lipschitz polyhedral surface ($\G = \partial\Omega$)
and the EFIE represents a boundary value problem for the time-harmonic Maxwell's
equations in the exterior domain ${\field{R}}^3 \backslash \Omega$,
then one can be more specific about the regularity
of the solution $\bu$ of {\rm (\ref{bie_var})}. In fact, the presence of
re-entrant edges is inevitable for the exterior of any polyhedral domain.
Therefore, using the results for Maxwell singularities
(see {\rm \cite{CostabelD_00_SEF}}), we conclude that for a closed polyhedral surface
$\G$ there exists $r \in (0,\frac 12)$ such that $\bu \in \bX^r$, and
$\bu \not\in \bX^{1/2}$.
\end{remark}

\begin{remark} \label{rem_singular}
The convergence rate of the $hp$-BEM for the EFIE is limited by the low
regularity of the solution to the corresponding Maxwell's equations,
especially for problems with screens which represent the least regular case.
To obtain sharp a priori error estimates, a refined approximation analysis
of singularities inherent to the solution of the EFIE is needed
(see {\rm \cite{BespalovH_NpB}}).
However, when the regularity of the solution to the EFIE is stated in terms of
Sobolev spaces of tangential vector fields on $\G$, the result of
Theorem~{\rm \ref{thm_converge}} is optimal with respect to both $h$ and~$p$.
\end{remark}

\section{Preliminaries} \label{sec_prelim}
\setcounter{equation}{0}

\subsection{Functional spaces, norms, and inner products} \label{sec_spaces}

In our previous paper \cite{BespalovH_Chp} we recalled definitions of
the full range of Sobolev spaces necessary for the convergence analysis of
the BEM for the EFIE (see \S{3.1} therein).
This included Sobolev spaces on a Lipschitz domain $\Omega\subset{\field{R}}^n$
and Sobolev spaces of scalar functions and tangential vector fields
on a piecewise smooth (open or closed) Lipschitz surface $\G \subset {\field{R}}^3$.
We have also defined basic differential operators on $\G$.
In the present paper we will use the same notation as in \cite{BespalovH_Chp}
for all differential operators, Sobolev spaces and their norms.
For convenience of the reader, let us repeat some
essential definitions, in particular, those for
Sobolev spaces of tangential vector fields.
Furthermore, we will introduce some more spaces and norms,
which are indispensable for the error analysis of the $hp$-BEM.

Let $\Omega$ be a Lipschitz domain in ${\field{R}}^n$.
We use a traditional notation for the Sobolev spaces $H^s(\Omega)$ ($s \ge -1$),
$H^s_0(\Omega)$ ($s \in (0,1]$), and $\tH^s(\Omega)$ ($s \in [-1,1]$)
with their standard norms (cf.,~\cite{LionsMagenes}).
In particular, for $s \in (0,1)$, the spaces $\tH^s(\Omega)$ are defined
by using the real K-method of interpolation, and 
for $s \in [-1,0)$, the spaces $H^s(\Omega)$, $\tH^s(\Omega)$ and
their norms are defined by duality
with $L^2(\Omega) = H^0(\Omega) = \tH^0(\Omega)$ as pivot space.
The norm and inner product in $L^2(\Omega)$ will be denoted as
$\|\cdot\|_{0,\Omega}$ and $\<\cdot,\cdot\>_{0,\Omega}$, respectively.

It is known that the standard norms $\|\cdot\|_{H^s(\Omega)}$
and $\|\cdot\|_{\tH^{-s}(\Omega)}$ are not scalable for $s \in (0,1]$
under affine transformations of $\Omega$ onto the reference domain (element).
However, we will need a scalable norm in the space $\tH^{-1}$ on a generic
interval $\ell \subset {\field{R}}^1$. 
Scalable families of norms in $H^s$ and in $\tH^{-s}$ for $s \in [-1,1]$
were introduced in \cite{ErvinH_06_ABS}. Following \cite{ErvinH_06_ABS} we
define
\[
  \|f\|^2_{H^1_h(\ell)} =
  (\hbox{meas}(\ell))^{-2} \|u\|^2_{0,\ell} + |u|^2_{H^1(\ell)}.
\]
Then the norm $\|\cdot\|_{\tH^{-1}_h(\ell)}$ is defined by duality:
\be \label{tH^{-1}_h-norm}
    \|f\|_{\tH^{-1}_h(\ell)} = \sup_{0 \not= \varphi \in H^{1}_h(\ell)}
    {|\<f,\varphi\>_{0,\ell}| \over{\|\varphi\|_{H^{1}_h(\ell)}}}.
\ee
Let $\ell$ be the image of the reference interval $\hat\ell$ under an affine transformation
$M$, i.e., $\ell = M(\hat\ell)$, and let $\hbox{meas}(\ell) \simeq h$.
Denote $\hat f = f \circ M$. Then the norms $\|\cdot\|_{H^{1}_h(\ell)}$ and
$\|\cdot\|_{\tH^{-1}_h(\ell)}$ are scalable (see \cite[Lemma 3.1]{ErvinH_06_ABS}),
i.e.,
\be \label{scale}
    \|f\|_{H^1_h(\ell)} \simeq h^{-1/2} \|\hat f\|_{H^1(\hat\ell)}\quad
    \hbox{and}\quad
    \|f\|_{\tH^{-1}_h(\ell)} \simeq h^{3/2} \|\hat f\|_{\tH^{-1}(\hat\ell)}
\ee
for any $\hat f \in H^1(\hat\ell)$ and $\hat f \in \tH^{-1}(\hat\ell)$, respectively,
and both equivalences are uniform for $h>0$.

An important fact related to the norm $\|\cdot\|_{\tH^{-1}_h}$
is that it enjoys the localisation property,
provided that the function has zero average on each sub-domain
(see \cite[Lemma 3.2]{ErvinH_06_ABS}). In particular, we will need the
following result.

\begin{lemma} \label{lm_split}
Assume that the boundary $\partial\Omega$ of the polygonal domain
$\Omega \subset {\field{R}}^2$ is partitioned into $N$ segments $\ell_j$ ($j=1,\ldots,N$).
Then, for all $f\in H^{-1}(\partial\Omega)$ with
$f|_{\ell_j} \in \tH^{-1}(\ell_j)$ and $\int_{\ell_j} f\,d\sigma = 0$ ($j=1,\ldots,N$),
there holds
\[
  \|f\|^2_{H^{-1}(\partial\Omega)}\le
  C \sum\limits_{j=1}^{N} \|f|_{\ell_j}\|^2_{\tH^{-1}_h(\ell_j)}
\]
with a positive constant $C$ independent of $f$ and $N$.
\end{lemma}

Now, let $\G$ be a piecewise smooth (open or closed) Lipschitz surface
in ${\field{R}}^3$. We will assume that $\G$ has plane faces $\G^{(i)}$
($i = 1,\ldots,\CI$; without loss of generality it is assumed that $\CI > 1$)
and straight edges $e_{ij} = \bar\G^{(i)} \cap \bar\G^{(j)} \not= \mbox{\o}$ ($i \not= j$).
If $\G$ is a closed surface, we will denote by
$\Omega$ the Lipschitz polyhedron bounded by $\G$, i.e., $\G=\partial\Omega$.
If $\G$ is an open surface, we additionally assume that $\G$ is orientable.
In this case, we first introduce a piecewise plane
closed Lipschitz surface $\tilde\G$ which contains $\G$, and then denote by
$\Omega$ the Lipschitz polyhedron bounded by $\tilde\G$, i.e., $\tilde\G=\partial\Omega$.
For each face $\G^{(i)} \subset \G$ there exists a constant unit normal vector $\bnu_i$,
which is an outer normal vector to $\Omega$.
These vectors are then blended into a unit normal vector
$\bnu$ defined almost everywhere on~$\G$.
For each pair of indices $i,j =  1,\ldots,\CI$ such that
$\bar\G^{(i)} \cap \bar\G^{(j)} = e_{ij}$
we consider unit vectors $\btau_{ij}$, $\btau_i^{(j)}$, and $\btau_j^{(i)}$
such that $\btau_{ij} \| e_{ij}$, $\btau_{i}^{(j)} = \btau_{ij} \times \bnu_i$,
and $\btau_{j}^{(i)} = \btau_{ij} \times \bnu_j$.
Since each $\G^{(i)}$ can be identified with a bounded subset in ${\field{R}}^2$,
the pair $(\btau_{i}^{(j)},\btau_{ij})$ is an orthonormal basis of the plane generated
by $\G^{(i)}$.

Let $\G$ be a closed surface. Then $\G=\partial\Omega$ is
locally the graph of a Lipschitz function. Since the Sobolev
spaces $H^s$ for $|s| \le 1$ are invariant under Lipschitz (i.e., $C^{0,1}$)
coordinate transformations, the spaces $H^s(\G)$ with $|s| \le 1$
are defined in the usual way via a partition of unity subordinate to a finite family
of local coordinate patches (see \cite{McLean_00_SES}).
Due to this definition, the properties of Sobolev spaces on Lipschitz domains
in ${\field{R}}^n$ carry over to Sobolev spaces on Lipschitz surfaces.
If $\G$ is an open surface, then the Sobolev spaces
$H^s(\G)$, $\tH^s(\G)$ for $|s| \le 1$ and $H^s_0(\G)$ for $0 < s \le 1$
are constructed in terms of the Sobolev spaces $H^s(\tilde\G)$ on
a closed Lipschitz surface $\tilde\G \supset \G$ (see \cite{McLean_00_SES}). 
Note that the spaces $H^s(\G^{(i)})$ and $\tH^s(\G^{(i)})$ on each face
$\G^{(i)}$ are well-defined for any $s \ge -1$.

We will denote by $\g_{\rm tr}$ the standard trace operator,
$\g_{\rm tr}(u) = u|_{\G}$, $u \in C^\infty(\bar\Omega)$.
For $s \in (0,1)$ (resp., $s > 1$), $\g_{\rm tr}$ has a unique extension
to a continuous operator $H^{s+1/2}(\Omega) \rightarrow H^{s}(\G)$
(resp., $H^{s+1/2}(\Omega) \rightarrow H^{1}(\G)$),
see \cite{Costabel_88_BIO,BuffaC_03_EFI}.
We will use the notation $C^\infty(\G) = \g_{\rm tr}(C^\infty(\bar\Omega))$.

Using the introduced Sobolev spaces of scalar functions, we define:
\[
  \bH^s(\Omega) = (H^s(\Omega))^3,\quad
  \bH^s(\G^{(i)}) = (H^s(\G^{(i)}))^2,\quad
  \tilde\bH^s(\G^{(i)}) = (\tH^s(\G^{(i)}))^2\quad
  (1 \le i \le \CI)
\]
for $s\ge -1$, and
\[
  \bH^s(\G) = (H^s(\G))^3\quad \hbox{for $s \in [-1,1]$}.
\]
The norms and inner products in all these spaces are defined component-wise and
usual conventions $\bH^0(\Omega) = \bL^2(\Omega)$,
$\bH^0(\G) = \bL^2(\G)$,
$\bH^0(\G^{(i)}) = \tilde\bH^0(\G^{(i)}) = \bL^2(\G^{(i)})$ hold.

Now let us introduce the Sobolev spaces of tangential vector fields
defined on $\G$ (see \cite{BuffaC_01_TFI,BuffaC_01_TII,BuffaCS_02_THL}).
We start with the space
\[
  \bL^2_t(\G) := \{\bu \in \bL^2(\G);\; \bu\cdot\bnu = 0\ \hbox{on}\ \G\},
\]
which will be identified with the space of two-dimensional, tangential,
square integrable vector fields. The norm and inner product in this space
will be denoted by $\|\cdot\|_{0,\G}$ and $\<\cdot,\cdot\>_{0,\G}$, respectively.
The similarity of this notation with the one for scalar functions
should not lead to any confusion. Then we define
(hereafter, $\bu_i$ denotes the restriction of $\bu$ to the face $\G^{(i)}$):
\beas
     & \bH_{\;-}^{s}(\G) :=
     \{\bu\in \bL^2_t(\G);\; \bu_i \in \bH^{s}(\G^{(i)}),\ \ 
     1\le i\le \CI\},\quad s \ge 0, &
     \cr\cr
     & \|\bu\|_{\bH_{\;-}^{s}(\G)} :=
     \bigg(\sum\limits_{i=1}^{\CI} \|\bu_i\|_{\bH^{s}(\G^{(i)})}^2\bigg)^{\frac 12}. &
\eeas

Let $\g_{\rm tr}$ be the trace operator (now acting on vector fields),
$\g_{\rm tr}(\bu) = \bu|_{\G}$,
$\g_{\rm tr}: \bH^{s+1/2}(\Omega) \rightarrow \bH^{s}(\G)$ for $s \in (0,1)$,
and let $\g_{\rm tr}^{-1}$ be one of its right inverses.
We will use the ``tangential components trace'' mapping
$\pi_\tau: (C^{\infty}(\bar\Omega))^3 \rightarrow \bL^2_t(\G)$ and
the ``tangential trace'' mapping
$\gamma_\tau: (C^{\infty}(\bar\Omega))^3 \rightarrow \bL^2_t(\G)$,
which are defined as $\bu \mapsto \bnu \times (\bu \times \bnu)|_{\G}$
and $\bu \times \bnu|_{\G}$, respectively.
We will also use the notation $\pi_\tau$ (resp., $\g_\tau$) for the composite
operator $\pi_\tau \circ \g_{\rm tr}^{-1}$ (resp., $\g_\tau \circ \g_{\rm tr}^{-1}$),
which acts on traces. Then we define the spaces
\[
  \bH^{1/2}_{\|}(\G) := \pi_\tau(\bH^{1/2}(\G)),\qquad
  \bH^{1/2}_{\perp}(\G) := \g_\tau(\bH^{1/2}(\G)),
\]
endowed with their operator norms
\beas
     \|\bu\|_{\bH^{1/2}_{\|}(\G)}
     & := &
     \inf_{\bphi \in \bH^{1/2}(\G)} \{\|\bphi\|_{\bH^{1/2}(\G)};\;
     \pi_\tau(\bphi) = \bu\},
     \cr\cr
     \|\bu\|_{\bH^{1/2}_{\perp}(\G)}
     & := &
     \inf_{\bphi \in \bH^{1/2}(\G)} \{\|\bphi\|_{\bH^{1/2}(\G)};\;
     \g_\tau(\bphi) = \bu\}.
\eeas
It has been shown in \cite{BuffaC_01_TFI} that the space $\bH^{1/2}_{\|}(\G)$
(resp., $\bH^{1/2}_{\perp}(\G)$) can be characterised as the space of tangential
vector fields belonging to $\bH_{\;-}^{1/2}(\G)$ and satisfying an appropriate
``weak continuity'' condition for the tangential (resp., normal) component
across each edge $e_{ij}$ of $\G$.

For $s > \frac 12$ we set
\beas
     \bH_{\|}^{s}(\G)
     & := &
     \{\bu \in \bH_{\;-}^{s}(\G);\;
     \bu_i \cdot \btau_{ij} = \bu_j \cdot \btau_{ij}\quad
     \hbox{at each}\ e_{ij}\},
     \cr\cr
     \bH_{\perp}^{s}(\G)
     & := &
     \{\bu \in \bH_{\;-}^{s}(\G);\;
     \bu_i \cdot \btau_{i}^{(j)} = \bu_j \cdot \btau_{j}^{(i)}\quad
     \hbox{at each}\ e_{ij}\}.
\eeas
For any $s > \frac 12$ the spaces $\bH_{\|}^{s}(\G)$ and $\bH_{\perp}^{s}(\G)$
are closed subspaces of $\bH_{\;-}^{s}(\G)$. Finally, for $s \in [0,\frac 12)$
we set
\[
  \bH_{\|}^{s}(\G) = \bH_{\perp}^{s}(\G) := \bH_{\;-}^{s}(\G).
\]

If $\G$ is an open surface, then we also need to define subspaces
of $\bH_{\|}^{s}(\G)$ and $\bH_{\perp}^{s}(\G)$ incorporating boundary
conditions on $\partial\G$ (for tangential and normal components, respectively).
In this case, for a given function $\bu$ on $\G$, we will denote by $\tilde\bu$
the extension of $\bu$ by zero onto a closed Lipschitz polyhedral surface
$\tilde\G \supset \G$. Then we define the spaces
\beas
     \tilde\bH_{\|}^{s}(\G)
     & := &
     \{\bu \in \bH_{\|}^{s}(\G);\;
     \tilde\bu \in \bH_{\|}^{s}(\tilde\G)\},\quad s \ge 0,
     \cr\cr
     \tilde\bH_{\perp}^{s}(\G)
     & := &
     \{\bu \in \bH_{\perp}^{s}(\G);\;
     \tilde\bu \in \bH_{\perp}^{s}(\tilde\G)\},\quad s \ge 0,
\eeas
which are furnished with the norms
\[
  \|\bu\|_{\tilde\bH^{s}_{\|}(\G)} :=
  \|\tilde\bu\|_{\bH^{s}_{\|}(\tilde\G)},\quad
  \|\bu\|_{\tilde\bH^{s}_{\perp}(\G)} :=
  \|\tilde\bu\|_{\bH^{s}_{\perp}(\tilde\G)},\quad s \ge 0.
\]
When considering open and closed surfaces at the same time we use the notation
$\tilde\bH_{\|}^{s}(\G)$, $\tilde\bH_{\perp}^{s}(\G)$, etc. also for closed
surfaces by assuming that $\tilde\bH_{\|}^{s}(\G) = \bH_{\|}^{s}(\G)$,
$\tilde\bH_{\perp}^{s}(\G) = \bH_{\perp}^{s}(\G)$, etc. in this case.
This in particular applies to the following definition of dual spaces.
For $s \in [-1,0)$, the spaces 
$\bH_{\|}^{s}(\G)$, $\tilde\bH_{\|}^{s}(\G)$, $\bH_{\perp}^{s}(\G)$, and
$\tilde\bH_{\perp}^{s}(\G)$ are defined as the dual
spaces of $\tilde\bH_{\|}^{-s}(\G)$, $\bH_{\|}^{-s}(\G)$,
$\tilde \bH_{\perp}^{-s}(\G)$, and $\bH_{\perp}^{-s}(\G)$,
respectively (with $\bL^2_t(\G)$ as pivot space).
They are equipped with their natural (dual) norms. Moreover, for any
$s \in (-\frac 12,\frac 12)$ there holds (cf. \cite{Grisvard_92_SBV})
\[
  \tilde\bH_{\|}^{s}(\G) = \bH_{\|}^{s}(\G) =
  \tilde\bH_{\perp}^{s}(\G) = \bH_{\perp}^{s}(\G).
\]

Using the above spaces of tangential vector fields, one can define basic
differential operators on $\G$. The tangential gradient,
$\gradg: H^1(\G) \rightarrow \bL^2_t(\G)$, and the tangential vector curl,
$\bcurlg: H^1(\G) \rightarrow \bL^2_t(\G)$, are defined in the usual way
by localisation to each face $\G^{(i)}$.
The adjoint operator of $-\gradg$ is the surface divergence denoted by $\divg$
(we refer to \cite{BuffaC_01_TFI,BuffaC_01_TII} for more details regarding
definitions and properties of differential operators on both closed and open surfaces).

Now we can introduce the spaces which appear when dealing with the EFIE on $\G$.
First, we set
\[
  \bH^s(\divg,\G) :=
  \{\bu \in \bH^s_{\|}(\G);\; \divg\,\bu \in H^s(\G)\},\quad s \in [-1/2,1/2]
\]
and
\[
  \bH^s_{\;-}(\divg,\G) :=
  \{\bu \in \bH^s_{\;-}(\G);\; \divg\,\bu \in H^s_{\;-}(\G)\},\quad s \ge 0
\]
for $\G$ being either a closed or an open surface.
Here, the space $H^{s}_{\;-}(\G)$ is defined similarly to the
space $\bH^{s}_{\;-}(\G)$ in a piecewise fashion:
\[
  H^s_{\;-}(\G) := \{u \in L^2(\G);\; u|_{\G^{(i)}} \in H^s(\G^{(i)}),\ 
                     i = 1,\ldots,\CI\},\quad s \ge 0,
\]
\[
  \|u\|^2_{H^s_{\;-}(\G)} := 
  \sum\limits_{i=1}^{\CI} \|u|_{\G^{(i)}}\|^2_{H^{s}(\G^{(i)})}.
\]
If $\G$ is an open surface, then we will also use the space
\[
  \tilde\bH^s(\divg,\G) :=
  \{\bu \in \tilde\bH^{s}_{\|}(\G);\;
  \divg\,\bu \in \tH^{s}(\G)\},\quad s \in [-1/2,0].
\]
The spaces $\bH^s(\divg,\G)$, $\bH^s_{\;-}(\divg,\G)$, and $\tilde\bH^s(\divg,\G)$
are equipped with their graph norms $\|\cdot\|_{\bH^s(\divg,\G)}$,
$\|\cdot\|_{\bH^s_{\;-}(\divg,\G)}$, and $\|\cdot\|_{\tilde\bH^s(\divg,\G)}$,
respectively. If $s=0$, then we will drop the superscript and for open surfaces also
the tilde in the above notation,
$\bH^0(\divg,\G) = \bH^0_{\;-}(\divg,\G) = \tilde\bH^0(\divg,\G) = \bH(\divg,\G)$.

On open surfaces, one needs the spaces incorporating homogeneous
boundary conditions for the trace of the normal component on $\partial\G$.
By $\tilde\bH^{s}_0(\divg,\G)$ with $s \in [-\frac 12,0]$ we denote the
subspace of elements $\bu \in \tilde\bH^{s}(\divg,\G)$
such that for all $v \in C^{\infty}(\G)$ there holds
\[
    \<\bu,\gradg v\> + \<\divg\,\bu,v\> = 0,
\]
where brackets $\<\cdot,\cdot\>$ denote the corresponding dualities.
Similarly as above, we will drop the superscript and the tilde if $s=0$.
We note that $\tilde\bH^{s}_0(\divg,\G)$ is a closed subspace of
$\tilde\bH^{s}(\divg,\G)$ for $s \in [-\frac 12,0]$.

The notation above related to open surfaces will be used also for
two-dimensional domains (in particular, for single faces of $\G$ and for
reference elements). When we need to join the notation for open and closed surfaces,
we will write
\[
  \bX = \bH^{-1/2}(\divg,\G),\quad
  \bX^{s} = \cases{
                   \bH^s(\divg,\G)         & \hbox{for $s \in [-\frac 12,0)$},
                   \cr
                   \noalign{\vspace{5pt}}
                   \bH(\divg,\G)          & \hbox{for $s = 0$},
                   \cr
                   \noalign{\vspace{5pt}}
                   \bH^s_{\;-}(\divg,\G)   & \hbox{for $s > 0$}
                   \cr
                  }
\]
if $\G$ is a closed surface, and
\[
  \bX = \tilde\bH^{-1/2}_0(\divg,\G),\quad
  \bX^{s} = \cases{
                   \tilde\bH^s_0(\divg,\G)                    & \hbox{for $s \in [-\frac 12,0)$},
                   \cr
                   \noalign{\vspace{5pt}}
                   \bH_0(\divg,\G)                            & \hbox{for $s = 0$},
                   \cr
                   \noalign{\vspace{5pt}}
                   \bH^s_{\;-}(\divg,\G) \cap \bH_0(\divg,\G) & \hbox{for $s > 0$}
                   \cr
                  }
\]
if $\G$ is an open surface. In all these cases the norms will be denoted
as $\|\cdot\|_{\bX}$ and $\|\cdot\|_{\bX^s}$.

Let $[\cdot,\cdot]_{\theta}$ ($\theta \in [0,1]$) denote the standard
interpolation (see \cite{LionsMagenes,BerghL_76_IS}). We quote the following interpolation
result from \cite{BuffaC_03_EFI} (see Theorem~4.12 therein).

\begin{lemma} \label{lm_interpolation}
There exists an $s_0 \in (0,\frac 12]$ such that for any $s \in [-\frac 12,s_0)$
there holds $[\bX,\bX^s]_{\theta} = \bX^{(1/2+s)\theta-1/2}$.
\end{lemma}

\subsection{Interpolation operators} \label{sec_int}

As it was mentioned in the introduction, our analysis of $hp$-approximations
essentially relies on the properties of the $\bH(\div)$-conforming projection based
interpolation operator. Let us briefly sketch the definition of this operator
on the reference element $K$ (see \cite{DemkowiczB_03_pIE} for details).
In this sub-section we use standard differential operators $\grad$, $\bcurl$ and $\div$
acting on 2D scalar functions and vector fields, respectively.

Given a vector field $\bu \in \bH^r(K) \cap \bH(\div,K)$ with $r > 0$,
the interpolant $\bu^p = \Pi^{\div}_p \bu \in \bCP^{\rm RT}_p(K)$ is defined as
the sum of three terms:
\[
  \bu^p = \bu_1 + \bu^p_2 + \bu^p_3.
\]
Here, $\bu_1$ is the lowest order interpolant
\[
  \bu_1 = \sum\limits_{\ell \subset \partial K}
  \Big(\int\limits_{\ell} \bu \cdot \hat\bn\,d\sigma\Big)\,\bphi_{\ell},
\]
where $\hat\bn$ denotes the outward normal unit vector to $\partial K$,
and $\bphi_{\ell}$ are the standard basis functions (associated with edges $\ell$)
for $\bCP_1^{\rm RT}(K)$.
For any edge $\ell \subset \partial K$ one has
\be \label{E-E_1}
    \int\limits_{\ell} (\bu - \bu_1) \cdot \hat\bn\,d\sigma = 0.
\ee
Hence, there exists a scalar function $\psi$, defined on $\partial K$,
such that
\be \label{psi}
    {\partial\psi\over{\partial \sigma}} =  (\bu - \bu_1) \cdot \hat\bn,\quad
    \psi = 0 \ \ \hbox{at all vertices}.
\ee
Then, for each edge $\ell$, the restriction $\psi|_{\ell}$ is projected
onto the set of polynomials $\CP^0_{p+1}(\ell)$
\be \label{psi_2^l}
    \psi_{2,\ell} \in \CP^0_{p+1}(\ell):\quad
    \<\psi|_{\ell} - \psi_{2,\ell},\varphi\>_{\tH^{1/2}(\ell)} = 0\quad
    \forall\,\varphi \in \CP^0_{p+1}(\ell).
\ee
Extending $\psi_{2,\ell}$ by zero from $\ell$ onto $\partial K$ (and keeping
its notation) and using the polynomial extension from the boundary,
we define $\psi_{2,p+1}^{\ell} \in \CP_{p+1}(K)$ such that
$\psi_{2,p+1}^{\ell}|_{\partial K} = \psi_{2,\ell}$. Then we set
\[
  \bu_2^p = \sum\limits_{\ell\, \subset\, \partial K}
            \bu^p_{2,\ell} \in \bCP_p^{\rm RT}(K),\ \ 
  \hbox{where \ } \bu^p_{2,\ell} = \bcurl\,\psi_{2,p+1}^{\ell}.
\]
The interior interpolant $\bu^p_3$ is a vector bubble function that solves
the constrained minimization problem
\[
  \ba{l}
  \bu^p_3 \in \bCP^{\rm RT,0}_p(K):
  \\[7pt]
  \|\div(\bu - (\bu_1 + \bu_2^p + \bu_3^p))\|_{0,K} \rightarrow \min,
  \\[7pt]
  \<\bu - (\bu_1 + \bu_2^p + \bu_3^p),\bcurl\,\phi\>_{0,K} = 0\quad
  \forall \phi \in \CP^{0}_{p+1}(K).
  \ea
\]

\begin{remark} \label{rem_tilde_1/2_ip}
The Sobolev space $\tH^{1/2}$ and the corresponding norm were defined in
{\rm \cite{BespalovH_Chp}}
using the real K-method of interpolation (see {\rm \cite{LionsMagenes}}).
However, the expression for the
$\tH^{1/2}$-inner product in {\rm (\ref{psi_2^l})} is based on another
(equivalent) norm in $\tH^{1/2}$. Without loss of generality, let us
assume that $\ell = I = (0,1)$ is the reference interval.
Then one has (see, e.g., {\rm \cite{DemkowiczB_03_pIE}})
\[
  \|\phi\|_{\tH^{1/2}(I)} \simeq 
  \|\grad\,\widetilde{\phi^\circ}\|_{0,T},
\]
where $T$ is the reference triangle, $\phi^\circ$ denotes the extension of the
function $\phi \in \tH^{1/2}(I)$ by zero onto $\partial T \supset I$,
and $\widetilde{\phi^\circ}$ is the harmonic lift of
$\phi^\circ \in H^{1/2}(\partial T)$. Then, applying the parallelogram
law and integrating by parts, we find the expression of the
corresponding inner product (cf. {\rm \cite{DemkowiczB_03_pIE}}):
\beas
     \<\phi,\psi\>_{\tH^{1/2}(I)}
     & = &
     \<\grad \widetilde{\phi^\circ},\grad \widetilde{\psi^\circ}\>_{0,T} =
     \Big\<
           \frac{\partial \widetilde{\phi^\circ}}{\partial \hat\bn},
           \widetilde{\psi^\circ}
     \Big\>_{0,\partial T} =
     \\[5pt]
     & = &
     \Big\<\frac{\partial \widetilde{\phi^\circ}}{\partial \hat\bn},\psi\Big\>_{0,I} =
     \Big\<\phi,\frac{\partial\widetilde{\psi^\circ}}{\partial \hat\bn}\Big\>_{0,I}\quad
     \forall \phi,\psi \in \tH^{1/2}(I).
\eeas
This expression can also be written as
\be \label{tilde_1/2_ip}
    \<\phi,\psi\>_{\tH^{1/2}(I)} =
    \Big\<\phi,\CDN_T(\psi)\Big\>_{0,I} = 
    \Big\<\CDN_T(\phi),\psi\Big\>_{0,I}\quad
    \forall \phi,\psi \in \tH^{1/2}(I),
\ee
where $\CDN_T: \tH^{1/2}(I) \rightarrow H^{-1/2}(I)$ denotes the Dirichlet-to-Neumann
operator associated with the triangle $T$.
Exactly this inner product is employed in {\rm (\ref{psi_2^l})}.
\end{remark}

For $r > 0$ the operator
$\Pi_p^{\div}:\; \bH^r(K) \cap \bH(\div,K) \rightarrow \bH(\div,K)$
is well defined and bounded, with corresponding operator norm being independent
of the polynomial degree $p$ (cf. \cite[Propositions~2]{DemkowiczB_03_pIE}).
Moreover, $\Pi_p^{\div}$ preserves polynomial vector fields
in $\bCP^{\rm RT}_p(K)$, and there holds the following estimate for
the interpolation error (see \cite[Theorem~5.1]{BespalovH_OEE})
\be \label{H(div)_p-estimate}
    \|\bu - \Pi_p^{\div}\,\bu\|_{\bH(\div,K)} \le
    C\, p^{-r}\, \|\bu\|_{\bH^r(\div,K)},
\ee
provided that
$\bu \in \bH^r(\div, K) := \{\bu \in \bH^r(K);\;\div\,\bu \in H^r(K)\}$
for $r > 0$.

We will need the following auxiliary result regarding the operator $\Pi_p^{\div}$.

\begin{lemma} \label{lm_aux1}
Let $\bu \in \bH^r(K) \cap \bH(\div,K)$ with $r > 0$, and let 
$\bu^p = \Pi^{\div}_p \bu \in \bCP^{\rm RT}_p(K)$. Then for any edge
$\ell \subset \partial K$ there holds
\be \label{aux1}
    \|(\bu - \bu^p) \cdot \hat\bn\|_{\tH^{-1}(\ell)} \le
    C\,p^{-1/2}\,\|(\bu - \bu^p) \cdot \hat\bn\|_{H^{-1/2}(\partial K)} \le
    C\,p^{-1/2}\,\|\bu - \bu^p\|_{H(\div,K)}.
\ee
\end{lemma}

\begin{proof}
Let us fix an edge $\ell \subset \partial K$.
Using the definition of the interpolant $\bu^p$ we have
\be \label{aux1_1}
    \|(\bu - \bu^p) \cdot \hat\bn\|_{\tH^{-1}(\ell)} =
    \|(\bu - \bu_1) \cdot \hat\bn - \bu_2^p \cdot \hat\bn\|_{\tH^{-1}(\ell)} =
    \Big\|
          \frac{\partial\psi}{\partial\sigma} -
          \frac{\partial\psi_{2,\ell}}{\partial\sigma}
    \Big\|_{\tH^{-1}(\ell)}.
\ee
Since the derivative with respect to the arc
length $\frac{\partial}{\partial\sigma}$ is a bounded operator
mapping $L^2(\ell)$ to $\tH^{-1}(\ell)$ (see, e.g., \cite[Lemma~3]{HahneS_96_SIE}),
we deduce from (\ref{aux1_1})
\be \label{aux1_2}
    \|(\bu - \bu^p) \cdot \hat\bn\|_{\tH^{-1}(\ell)} \le
    C\,\|\psi|_{\ell} - \psi_{2,\ell}\|_{0,\ell}.
\ee
Let $\tilde T$ be an equilateral triangle having $\ell$ as one of its edges
(if $K = T$ then $\tilde T = T$). Denoting $\phi := \psi|_{\ell} - \psi_{2,\ell}$
and using the expression for the $\tH^{1/2}(\ell)$-inner product as in (\ref{tilde_1/2_ip}),
we can rewrite the orthogonality relation in (\ref{psi_2^l}) as
\be \label{aux1_3}
    \<\phi,\phi_p\>_{\tH^{1/2}(\ell)} =
    \Big\<\phi,\CDN_{\tilde T}(\phi_p)\Big\>_{0,\ell} = 0\qquad
    \forall\,\phi_p \in \CP_{p+1}^0(\ell).
\ee
Now, let us consider an auxiliary mixed boundary value problem on $\tilde T$:
find $\Phi \in H^1(\tilde T)$ such that
\[
  -\Delta \Phi = 0\ \ \hbox{in $\tilde T$},\quad
  \frac{\partial\Phi}{\partial \hat\bn} = \phi\ \ \hbox{on $\ell$},\quad
  \Phi = 0\ \ \hbox{on $\partial \tilde T \backslash \ell$}.
\]
Recalling that $\phi \in \tH^{1/2}(\ell) \subset H^{1/2}(\ell)$ and
using the regularity theory for elliptic problems in non-smooth domains
(see, e.g., \cite{Dauge_88_EBV,Grisvard_92_SBV}), we conclude that
$\Phi \in H^2(\tilde T)$ and $\|\Phi\|_{H^2(\tilde T)} \le C\,\|\phi\|_{H^{1/2}(\ell)}$.
Therefore, applying the trace theorem for a single edge $\ell \subset \tilde T$
(see \cite[Theorem~1.4.2]{Grisvard_92_SBV}) we prove that
$\Phi|_{\ell} \in H^{3/2}(\ell)$ and
\be \label{aux1_4}
    \|\Phi|_{\ell}\|_{H^{3/2}(\ell)} \le
    C\,\|\Phi\|_{H^2(\tilde T)} \le
    C\,\|\phi\|_{H^{1/2}(\ell)}.
\ee
On the other hand,
$\phi = {\partial\Phi}/{\partial \hat\bn} = \CDN_{\tilde T}(\Phi|_{\ell})$ and
using (\ref{aux1_3}) we obtain for any $\phi_p \in \CP_{p}^0(\ell)$
\beas
     \|\phi\|^2_{0,\ell}
     & = &
     \<\phi,\phi\>_{0,\ell} =
     \Big\<\phi,\CDN_{\tilde T}(\Phi|_{\ell})\Big\>_{0,\ell} =
     \Big\<\phi,\CDN_{\tilde T}(\Phi|_{\ell} - \phi_p)\Big\>_{0,\ell}
     \\[5pt]
     & \le &
     \|\phi\|_{0,\ell}\,\,\Big\|\CDN_{\tilde T}(\Phi|_{\ell} - \phi_p)\Big\|_{0,\ell}.
\eeas
Hence, using the fact that the operator $\CDN_{\tilde T}$ is continuous
as a mapping $H^1_0(\ell) \rightarrow L^2(\ell)$, we find
\be \label{aux1_5}
    \|\phi\|_{0,\ell} \le
    C\,|\Phi|_{\ell} - \phi_p|_{H^1(\ell)}\qquad
    \forall\,\phi_p \in \CP_{p}^0(\ell).
\ee
To estimate the semi-norm in (\ref{aux1_5}) we apply the standard
$p$-approximation result in 1D (see, e.g., \cite[Lemma~3.2]{BabuskaS_87_OCR}):
there exists $\phi_p \in \CP_{p}^0(\ell)$ such that
\be \label{aux1_6}
    \|\Phi|_{\ell} - \phi_p\|_{H^1(\ell)} \le
    C\,p^{-1/2}\,\|\Phi|_{\ell}\|_{H^{3/2}(\ell)}.
\ee
Putting together inequalities (\ref{aux1_4})--(\ref{aux1_6}) and recalling
our notation for the function $\phi$, we obtain
\be \label{aux1_7}
    \|\psi|_{\ell} - \psi_{2,\ell}\|_{0,\ell} \le
    C\,p^{-1/2}\,\|\psi|_{\ell} - \psi_{2,\ell}\|_{H^{1/2}(\ell)}.
\ee
Since
\[
  (\bu - \bu^p) \cdot \hat\bn|_{\partial K} =
  (\bu - \bu_1 - \bu_2^p) \cdot \hat\bn|_{\partial K} =
  \frac{\partial}{\partial\sigma}
  \Big(\psi - \sum_{\ell \subset \partial K} \psi_{2,\ell}\Big) \in
  H^{-1/2}_{*}(\partial K)
\]
with $H^{-1/2}_{*}(\partial K) := \{u \in H^{-1/2}(\partial K);\;\<u,1\>_{0,\partial K} = 0\}$,
and the tangential derivative defines an isomorphism
$\frac{\partial}{\partial\sigma}:\,
 H^{1/2}(\partial K)/{\field{R}} \rightarrow H^{-1/2}_{*}(\partial K)$
(see \cite[Lemma~2]{DemkowiczB_03_pIE}), we prove that
\be \label{aux1_8}
    \|\psi|_{\ell} - \psi_{2,\ell}\|_{H^{1/2}(\ell)} \le
    C\,\Big\|\psi - \sum_{\ell \subset \partial K} \psi_{2,\ell}\Big\|_{H^{1/2}(\partial K)} \le
    C\,\|(\bu - \bu^p) \cdot \hat\bn\|_{H^{-1/2}(\partial K)}.
\ee
The first inequality in (\ref{aux1}) then immediately follows from
estimates (\ref{aux1_2}), (\ref{aux1_7}), and (\ref{aux1_8}). The second inequality
in (\ref{aux1}) is true due to the continuity of the normal trace operator
$\bv \rightarrow \bv \cdot \hat\bn|_{\partial K}$ as a mapping
$\bH(\div,K) \rightarrow H^{-1/2}(\partial K)$.
\end{proof}

Now, let us consider our Lipschitz surface $\G$ discretised by the quasi-uniform mesh.
Using the Piola transform $\CM_j$ for each element $\G_j$ and applying the
$p$-interpolation operator $\Pi_p^{\div}$ on the reference elements,
one can define the ``global'' $\bH(\div)$-conforming $hp$-interpolation operator
$\Pi_{hp}^{\div}:\, \bH^r_{\;-}(\G) \cap \bH(\divg,\G) \rightarrow \bX_{hp}$
($r > 0$) such that for $\bu^{hp} := \Pi_{hp}^{\div} \bu$ there holds
\[
  \CM_j^{-1}(\bu^{hp}|_{\G_j}) = \Pi_p^{\div}\Big(\CM_j^{-1}(\bu|_{\G_j})\Big).
\]
Since the projection-based interpolation operator $\Pi^{\div}_p$
preserves polynomial vector fields and provides $\bH(\div)$-conforming
approximations, the $p$-interpolation error estimate (\ref{H(div)_p-estimate})
on the reference element extends to the corresponding $hp$-estimate
in a standard way by using the Bramble-Hilbert argument and scaling.
This result is formulated in the following theorem.

\begin{theorem} \label{thm_H(div)_hp-estimate}
Let $\bu \in \bH^r_{\;-}(\divg,\G)$, $r > 0$. Then there exists a positive constant $C$
independent of $h$, $p$, and $\bu$ such that
\[
  \|\bu - \Pi_{hp}^{\div}\bu\|_{\bH(\divg,\G)} \le
  C\,h^{\min\,\{r,p\}}\,p^{-r}\,\|\bu\|_{\bH^r_{\;-}(\divg,\G)}.
\]
\end{theorem}

\section{Approximating properties of $\bX_{hp}$ in $\bX$} \label{sec_approx}
\setcounter{equation}{0}

The main purpose of this section is to prove that the orthogonal projection
$P_{hp}:\,\bX \rightarrow \bX_{hp}$ with respect to the norm in $\bX$ satisfies
an optimal error estimate (when the error is measured in the norm of $\bX$).
The standard approach to tackle such kind of problems is to use the duality
argument. Following this approach on $\G$ in our setting would require
that for $s \in (0,\frac 12]$ the
space $\bH^{-s}(\divg,\G)$ (or, $\tilde\bH^{-s}(\divg,\G)$ if $\G$ is an open
surface) is the dual space of $\bH^{s}(\divg,\G)$ with respect to the
$\bH(\divg,\G)$-inner product. However, this fact is true only for smooth surfaces.
On the other hand, the duality argument can be applied face by face. This idea
was suggested by Buffa \& Christiansen in \cite{BuffaC_03_EFI} and was successfully
exploited by these authors in the context of the $h$-version of the BEM for the
EFIE. We will demonstrate below that, when employing the projection based
$\bH(\div)$-conforming interpolation operator, this approach provides an optimal
$hp$-approximation result in the energy norm of the EFIE.

Let $\tG$ be a single face of $\G$. To simplify the presentation (in particular,
to avoid imposing boundary conditions on the edges of $\partial\tG$ which also
belong to $\partial\G$ in the case that $\partial\G\not=\emptyset$), we will assume that $\G$
is a (closed) Lipschitz polyhedral surface.
All arguments extend to the case of an open piecewise plane
Lipschitz surface in a straightforward way (cf. \cite{BuffaC_03_EFI}).
For the sake of simplicity of notation we will also omit the subscript $\tG$
for differential operators over this face, e.g., we will write $\div$ for $\div_{\tG}$.

Let $\bH^r(\div, \tG) := \{\bu \in \bH^r(\tG);\;\div\,\bu \in H^r(\tG)\}$
for $r \ge 0$.  We will also denote by $\bX_{hp}(\tG)$ the restriction
of $\bX_{hp}$ to $\tG$. Then, for $r > 0$ we introduce
the operator $\CQ_{hp}:\,\bH^r(\div,\tG) \rightarrow \bX_{hp}(\tG)$ as follows:
given $\bu \in \bH^r(\div,\tG)$, we define $\CQ_{hp}\bu \in \bX_{hp}(\tG)$ such that
\be \label{def_Q_hp}
    \ba{rl}
    (\bu - \CQ_{hp}\bu,\bv)_{\bH(\div,\tG)} = 0\quad &
    \forall\,\bv \in \bX_{hp}(\tG) \cap \bH_0(\div,\tG),
    \\[10pt]
    \CQ_{hp}\bu \cdot \tbn = \Pi_{hp}^{\div}\bu \cdot \tbn\quad &
    \hbox{on $\partial\tG$}.
    \ea
\ee
Here, $(\cdot,\cdot)_{\bH(\div,\tG)}$ denotes the $\bH(\div,\tG)$-inner product,
and $\tbn$ is the unit outward normal vector to $\partial\tG$.
It follows immediately from (\ref{def_Q_hp}) that
\be \label{Q_hp_est1}
    \|\bu - \CQ_{hp}\bu\|_{\bH(\div,\tG)} \le
    \|\bu - \Pi_{hp}^{\div}\bu\|_{\bH(\div,\tG)}.
\ee

\begin{lemma} \label{lm_Q_hp}
For $r>0$ let $\bu \in \bH^r(\div,\tG)$.
Then there holds
\be \label{Q_hp_est2}
    \|\bu - \CQ_{hp}\bu\|_{\tilde\bH^{-1/2}(\div,\tG)} \le
    C \left(\hbox{$\frac hp$}\right)^{1/2}
    \|\bu - \Pi_{hp}^{\div}\bu\|_{\bH(\div,\tG)},
\ee
where the constant $C$ is independent of $h$ and $p$.
\end{lemma}

\begin{proof}
We follow the technique used by Buffa and Christiansen in the proof of Proposition~4.6
in~\cite{BuffaC_03_EFI} but rely on the properties of the $\bH(\div)$-conforming projection based
interpolation operator $\Pi_{hp}^{\div}$. For given $r>0$ and $\bu \in \bH^r(\div,\tG)$
let us consider the following problem: find $\bu_0 \in \bH(\div,\tG)$ such that
\beas
     (\bu - \bu_0,\bv)_{\bH(\div,\tG)} = 0
     & \quad &
     \forall\,\bv \in \bH_0(\div,\tG),
     \\[5pt]
     \bu_0 \cdot \tbn = \Pi_{hp}^{\div}\bu \cdot \tbn
     & \quad &
     \hbox{on $\partial\tG$}.
\eeas
Then, as an immediate consequence of Lemma~4.8 in \cite{BuffaC_03_EFI}, there holds
\be \label{Q_hp_est2_1}
    \|\bu - \bu_0\|_{\tilde\bH^{s+1/2}(\div,\tG)} \le
    C\,\|(\bu - \Pi_{hp}^{\div}\bu) \cdot \tbn\|_{H^{s}(\partial\tG)},\quad
    \hbox{$s = -1,\; -\frac 12$}.
\ee
By the triangle inequality one has
\be \label{Q_hp_est2_2}
    \|\bu - \CQ_{hp}\bu\|_{\tilde\bH^{-1/2}(\div,\tG)} \le
    \|\bu - \bu_0\|_{\tilde\bH^{-1/2}(\div,\tG)} +
    \|\bu_0 - \CQ_{hp}\bu\|_{\tilde\bH^{-1/2}(\div,\tG)}.
\ee
Let us estimate each term on the right-hand side of (\ref{Q_hp_est2_2}).
For the first term we have by using (\ref{Q_hp_est2_1}) with $s = -1$
\be \label{Q_hp_est2_3}
    \|\bu - \bu_0\|_{\tilde\bH^{-1/2}(\div,\tG)} \le
    C\,\|(\bu - \Pi_{hp}^{\div}\bu) \cdot \tbn\|_{H^{-1}(\partial\tG)}.
\ee
Denote $\bu^{hp} := \Pi_{hp}^{\div}\bu$. Since
$\int_{\ell_h} (\bu - \bu^{hp}) \cdot \tbn\, d\sigma = 0$
for any mesh edge $\ell_h \subset \partial\tG$, we can apply the localisation
result of Lemma~\ref{lm_split} to estimate
\be \label{Q_hp_est2_4}
    \|(\bu - \bu^{hp}) \cdot \tbn\|^2_{H^{-1}(\partial\tG)}\le
    C \sum\limits_{\ell_h \subset \partial\tG}
    \|(\bu - \bu^{hp}) \cdot \tbn\|^2_{\tH^{-1}_h(\ell_h)}.
\ee
Let us fix an edge $\ell_h \subset \partial\tG$. Then there exists an element
$K_h \subset \tG$ such that $\bar\ell_h = \partial\tG \cap \partial K_h$.
This element $K_h$ is the image of the reference element $K$ under an affine
mapping $M$, and let $\CM$ be the corresponding Piola transform.
Then using (\ref{tH^{-1}_h-norm}), (\ref{scale}) and the standard property
of the Piola transform (see \cite[Lemma~1.5]{BrezziF_91_MHF}), we have
\beas
     \|(\bu - \bu^{hp}) \cdot \tbn\|_{\tH^{-1}_h(\ell_h)}
     & = &
     \sup_{0 \not= \varphi \in H^{1}_h(\ell_h)}
     {\<(\bu - \bu^{hp})\cdot\tbn,\varphi\>_{0,\ell_h} \over{\|\varphi\|_{H^{1}_h(\ell_h)}}}
     \\[3pt]
     & \simeq &
     \sup_{0 \not= \hat\varphi \in H^{1}(\hat\ell)}
     {\<(\hat\bu - \Pi_p^{\div}\hat\bu)\cdot\hat\bn,\hat\varphi\>_{0,\hat\ell}
     \over{h^{-1/2}\,\|\hat\varphi\|_{H^{1}(\hat\ell)}}}
     \\[3pt]
     & = &
     C\,h^{1/2}\,\|(\hat\bu - \Pi_p^{\div}\hat\bu) \cdot \hat\bn\|_{\tH^{-1}(\hat\ell)},
\eeas
where $\hat\ell = M^{-1}(\ell_h) \subset \partial K$,
$\hat\bu = \CM^{-1}(\bu|_{K_h})$,
$\hat\varphi = \varphi \circ M$, and $\hat\bn$ denotes the unit outward normal
vector to $\partial K$. Hence, applying Lemma~\ref{lm_aux1} and using
standard scaling properties of the Piola transform
(see \cite[Lemma~1.6]{BrezziF_91_MHF}), we estimate
\be \label{Q_hp_est2_5}
    \|(\bu - \bu^{hp}) \cdot \tbn\|_{\tH^{-1}_h(\ell_h)} \le
    C \left(\hbox{$\frac hp$}\right)^{1/2}
    \|\hat\bu - \Pi_p^{\div}\hat\bu\|_{\bH(\div,K)} \le
    C \left(\hbox{$\frac hp$}\right)^{1/2}
    \|\bu - \bu^{hp}\|_{\bH(\div,K_h)}.
\ee
Combining inequalities (\ref{Q_hp_est2_5}) over all mesh edges
$\ell_h \subset \partial\tG$ and recalling the notation for $\bu^{hp}$, we deduce from
(\ref{Q_hp_est2_3}) and (\ref{Q_hp_est2_4}):
\be \label{Q_hp_est2_6}
    \|\bu - \bu_0\|_{\tilde\bH^{-1/2}(\div,\tG)} \le
    C \left(\hbox{$\frac hp$}\right)^{1/2}
    \|\bu - \Pi_{hp}^{\div}\bu\|_{\bH(\div,\tG)}.
\ee
Now we focus on the second term on the right-hand side of (\ref{Q_hp_est2_2})
and prove that
\be \label{Q_hp_est2_7}
    \|\bu_0 - \CQ_{hp}\bu\|_{\tilde\bH^{-1/2}(\div,\tG)} \le
    C \left(\hbox{$\frac hp$}\right)^{1/2}
    \|\bu - \Pi_{hp}^{\div}\bu\|_{\bH(\div,\tG)}.
\ee
Denote $\bX_{\tG} := \tilde\bH^{-1/2}_0(\div,\tG)$ and
$\bH^{1/2}_0(\div,\tG) := \bH^{1/2}(\div,\tG) \cap \bH_0(\div,\tG)$.
Let $\bX'_{\tG}$ be the dual space of $\bX_{\tG}$.
Then, by Lemma~4.7 in \cite{BuffaC_03_EFI}, the operator
$I - \grad(\div):\, \bH_0^{1/2}(\div,\tG) \rightarrow \bX'_{\tG}$ is an
isomorphism. Moreover, it is easy to see that $\CQ_{hp}\bu_0 = \CQ_{hp}\bu$ and
$\bu_0 - \CQ_{hp}\bu_0 \in \bH_0(\div,\tG) \subset \bX_{\tG}$.
Therefore, we can estimate as follows:
\beas
     \|\bu_0 - \CQ_{hp}\bu\|_{\tilde\bH^{-1/2}(\div,\tG)}
     & = &
     \|\bu_0 - \CQ_{hp}\bu_0\|_{\tilde\bH^{-1/2}(\div,\tG)}
     \\[3pt]
     & \le &
     C\, \sup_{\bzero \not= \bv \in \bX'_{\tG}}
     {\<\bv, \bu_0 - \CQ_{hp}\bu_0\>_{\bX'_{\tG},\bX_{\tG}}
     \over{\|\bv\|_{\bX'_{\tG}}}}
     \\[3pt]
     & \le &
     C\,\sup_{\bzero \not= \bw \in \bH_0^{1/2}(\div,\tG)}
     {\<\bw - \grad(\div\,\bw), \bu_0 - \CQ_{hp}\bu_0\>_{\bX'_{\tG},\bX_{\tG}}
     \over{\|\bw\|_{\bH^{1/2}(\div,\tG)}}}
     \\[3pt]
     & = &
     C\,\sup_{\bzero \not= \bw \in \bH_0^{1/2}(\div,\tG)}
     {(\bu_0 - \CQ_{hp}\bu_0, \bw)_{\bH(\div,\tG)}
     \over{\|\bw\|_{\bH^{1/2}(\div,\tG)}}}
     \\[3pt]
     & = &
     C\,\sup_{\bzero \not= \bw \in \bH_0^{1/2}(\div,\tG)}
     {(\bu_0 - \CQ_{hp}\bu_0, \bw - \Pi_{hp}^{\div}\bw)_{\bH(\div,\tG)}
     \over{\|\bw\|_{\bH^{1/2}(\div,\tG)}}};
\eeas
for the last step we used the definition of $\CQ_{hp}$
(see (\ref{def_Q_hp}) with $\bv = \Pi_{hp}^{\div}\bw$).
Hence, using the Cauchy-Schwarz inequality and the interpolation error estimate of
Theorem~\ref{thm_H(div)_hp-estimate} (with $\G = \tG$ and $r = \frac 12$)
we prove that
\be \label{Q_hp_est2_8}
    \|\bu_0 - \CQ_{hp}\bu\|_{\tilde\bH^{-1/2}(\div,\tG)} \le
    C \left(\hbox{$\frac hp$}\right)^{1/2}
    \|\bu_0 - \CQ_{hp}\bu\|_{\bH(\div,\tG)}.
\ee
The norm on the right-hand side of (\ref{Q_hp_est2_8}) is estimated by applying
the triangle inequality, (\ref{Q_hp_est2_1}) with $s = -\frac 12$, (\ref{Q_hp_est1}),
and the continuity property of the normal trace operator:
\bea
    \|\bu_0 - \CQ_{hp}\bu\|_{\bH(\div,\tG)}
    & \le &
    \|\bu_0 - \bu\|_{\bH(\div,\tG)} +
    \|\bu - \CQ_{hp}\bu\|_{\bH(\div,\tG)}
    \nonumber
    \\[2pt]
    & \le &
    C\,\|(\bu - \Pi_{hp}^{\div}\bu) \cdot \tbn\|_{H^{-1/2}(\partial\tG)} +
    \|\bu - \Pi_{hp}^{\div}\bu\|_{\bH(\div,\tG)}
    \nonumber
    \\[2pt]
    & \le &
    C\,\|\bu - \Pi_{hp}^{\div}\bu\|_{\bH(\div,\tG)}.
    \label{Q_hp_est2_9}
\eea
The desired inequality in (\ref{Q_hp_est2_7}) then follows from
(\ref{Q_hp_est2_8}) and (\ref{Q_hp_est2_9}).

To obtain (\ref{Q_hp_est2}) it remains to collect (\ref{Q_hp_est2_6})
and (\ref{Q_hp_est2_7}) in (\ref{Q_hp_est2_2}). This finishes the proof.
\end{proof}

Now we can prove the main result of this section.

\begin{theorem} \label{thm_main_approx}
Let $P_{hp}:\, \bX \rightarrow \bX_{hp}$ be the orthogonal projection
with respect to the norm in $\bX$. If $\bu \in \bX^r$ with
$r > -\frac 12$, then
\be \label{main_approx}
    \|\bu - P_{hp}\bu\|_{\bX} \le
    C\,h^{1/2+\min\,\{r,p\}}\,p^{-(r+1/2)}\,\|\bu\|_{\bX^r}
\ee
with a positive constant $C$ independent of $h$, $p$, and $\bu$.
\end{theorem}

\begin{proof}
First, let us assume that $r > 0$. We consider the discrete vector field
$\bv \in \bX_{hp}$ such that $\bv|_{\G^{(i)}} = \CQ^{(i)}_{hp}(\bu|_{\G^{(i)}})$
for each face $\G^{(i)}$ of $\G$ (here, $\CQ^{(i)}_{hp}$ denotes the operator defined as in
(\ref{def_Q_hp}) with respect to the face $\G^{(i)}$). Due to the localisation properties
of the norms in $\bX$ and in $\bH(\divg,\G)$, we have by Lemma~\ref{lm_Q_hp}
\beas
     \|\bu - P_{hp}\bu\|_{\bX}
     & \le &
     \|\bu - \bv\|_{\bX} \le
     C\,\sum\limits_{i=1}^{\CI}
     \|
       \bu|_{\G^{(i)}} -
       \CQ^{(i)}_{hp}(\bu|_{\G^{(i)}})
     \|_{\tilde\bH^{-1/2}(\div_{\G^{(i)}},\G^{(i)})}
     \\[3pt]
     & \le &
     C \left(\hbox{$\frac hp$}\right)^{1/2}
     \|\bu - \Pi_{hp}^{\div}\bu\|_{\bH(\divg,\G)},
\eeas
and inequality (\ref{main_approx}) then follows from the error
estimate of Theorem~\ref{thm_H(div)_hp-estimate}.

Now, let $r \in (-\frac 12, 0]$. Assume that $\bu \in \bX^s$ with
some $s \in (0,s_0)$, where $s_0 \in (0,\frac 12]$ is the same as
in Lemma~\ref{lm_interpolation}. Then, using the first part of the proof, one has
\[
  \|\bu - P_{hp}\bu\|_{\bX} \le
  C \left(\hbox{$\frac hp$}\right)^{1/2+s}
  \|\bu\|_{\bX^s}.
\]
On the other hand, it is trivial that
\[
  \|\bu - P_{hp}\bu\|_{\bX} \le \|\bu\|_{\bX}.
\]
Therefore, applying the interpolation argument which relies on Lemma~\ref{lm_interpolation},
we prove
\[
  \|\bu - P_{hp}\bu\|_{\bX} \le
  C \left(\hbox{$\frac hp$}\right)^{1/2+r}
  \|\bu\|_{\bX^r}\quad
  \forall\,\bu \in \bX^s.
\]
This estimate yields (\ref{main_approx}) due to the density of regular functions in
$\bX^r$, and the proof is finished.
\end{proof}

\bigskip\medskip

\noindent{\bf Acknowledgement.}
A significant part of this work has been done while A.B. was visiting
the Facultad de Matem\'aticas, Pontificia Universidad Cat\'olica de Chile
(Santiago, Chile). Their hospitality is gratefully acknowledged.



\end{document}